\documentclass[12pt,thmsa]{article}
\usepackage{amsfonts}

\usepackage{sw20jart}
\usepackage{epsfig}

\input{tcilatex}
\begin{document}

\author{R. F. Picken \\
Departamento de Matem\'{a}tica and\\
Centro de Matem\'{a}tica Aplicada\\
Instituto Superior T\'{e}cnico, Avenida\\
Rovisco Pais\\
1049-001 Lisboa, Portugal\\
and\\
P. A. Semi\~{a}o \\
\'{A}rea Departamental de Matem\'{a}tica and\\
Centro de Matem\'{a}tica Aplicada\\
Universidade do Algarve,\\
Unidade de Ci\^{e}ncias Exactas e Humanas,\\
8000 Gambelas, Portugal}
\title{{\Large \textbf{TQFT - a new direction in algebraic topology}}\thanks{%
Contribution to the proceedings of the workshop ``New Developments in
Algebraic Topology'', Faro, Portugal, July 13-14, 1998.}}
\maketitle

\begin{abstract}
We give an introduction for the non-expert to TQFT (Topological Quantum
Field Theory), focussing especially on its role in algebraic topology. We
compare the Atiyah axioms for TQFT with the Eilenberg Steenrod axioms for
homology, give a few simple examples of TQFTs, and discuss some other
approaches that have been taken to defining TQFT. We then propose a new
formulation of TQFT, which is closer in spirit to the way conventional
functors of algebraic topology, like homology, are presented. In this
approach the fundamental operation of gluing is incorporated through the
notion of a \emph{gluing morphism}, which we define. It allows not only the
gluing together of two separate objects, but also the self-gluing of a
single object to be treated in the same fashion. As an example of our
approach we reformulate and generalize a class of examples due to Quinn
based on the Euler characteristic.
\vspace{0.5cm}

\begin{flushright}
IST/DM/36/99\\
math.QA/9912085\\
\end{flushright}

\vspace{0.5cm}
\end{abstract}

\section{Introduction}

Topological Quantum Field Theory, or TQFT for short, is a notion which
originally arose from ideas of quantum physics. Since then the notion has
developed considerably in a number of directions, and in particular has had
a pervasive influence in mathematics. In this article we will be defending
the point of view that TQFT is a new type of functor of algebraic topology,
analogous to homology or homotopy.

We recall that a functor is a map between two categories, like the functor $%
H_{n}$ from the category of topological spaces to the category of abelian
groups, which assigns to each topological space $X$ its $n$-th homology
group $H_{n}(X)$, and \emph{furthermore} assigns to each continuous map $X%
\stackrel{f}{\rightarrow }Y$ a group morphism $H_{n}(X)\stackrel{H_{n}(f)}{%
\rightarrow }H_{n}(Y)$, with the assignments obeying certain natural
conditions. This property of not only mapping mathematical objects, but also
the morphisms between them, implies a very profound relationship between the
two categories, a bridge between two often very different regions of the
mathematical landscape.

For those without a physics background, a few words are also in order about
the underlying notions of quantum physics. Fields are simply (a suitable
class of) functions $\phi :\Sigma \rightarrow X$ defined on the $d$
-dimensional space manifold $\Sigma $, or $\Phi :M\rightarrow X$ defined on
the $d+1$-dimensional space-time manifold $M$, into another space $X$, which
may be, for instance, $\mathbf{R}^{n}$, or a Lie group $G$. A field $\Phi $
describes an evolution of fields $\phi $, at least when $M=\Sigma \times I$,
where $I$ is an interval (of time), and is interpreted as some kind of
multifingered evolution, when $M$ is not of this form. The evolution is
classical when $\Phi $ is a critical point of a certain functional $S(\Phi )$%
, called the action. The corresponding quantum evolution is a superposition
of evolutions known as the path integral, $\int d\mu (\Phi )\exp (-iS(\Phi
)/\hbar )$, where $d\mu (\Phi )$ is a rather elusive ``measure'' on the
space of all evolutions, and $\hbar $ is the quantum parameter. Roughly
speaking, as this parameter tends to zero, the path integral reduces to the
classical evolution in the stationary phase approximation. When $S(\Phi )$
and $d\mu (\Phi )$ happen to be invariant under diffeomorphisms, we have a
so-called \emph{topological} quantum field theory, giving rise to
topological invariants of $M$.

The most famous example of a TQFT was Witten's use of the Chern-Simons
action to obtain a topological invariant for $3$-manifolds (as well as,
simultaneously, an invariant of embedded $1$-dimensional submanifolds, i.e.,
knots and links) \cite{Witten}. Numerous other constructions, both heuristic
and rigorous, followed, including a class of ``state-sum models'', involving
piecewise-linear (PL) manifolds. \cite{Turaev,Turaev-Viro,ResTuraev}. 
For PL manifolds of dimension $4$, a state-sum model due to
Crane and Yetter \cite{CraneYetter} gave rise to a combinatorial formula for
the signature. TQFT has had most impact in $3$ and $4$-dimensional topology,
where the classical invariants are weak. For reviews and further background
see \cite{Lawrence,Sawin}.

\section{Atiyah's axioms}

The common features underlying a number of different TQFT constructions were
formalized by Atiyah in a set of axioms for TQFT \cite{Atiyah}. These are
modeled on similar axioms for conformal field theory, due to Segal \cite
{Segal}. We will give a brief description of them. According to Atiyah, a $%
(d+1)$-dimensional TQFT is an assignment $\Sigma \mapsto V_{\Sigma }$ and $%
M\mapsto Z_{M}$, which assigns to every $d$-dimensional oriented manifold $%
\Sigma $ a finite-dimensional vector space $V_{\Sigma }$ (over some fixed
field $\Bbb{K}$), and to every $(d+1)$-dimensional oriented manifold with
boundary $M$, an element of $V_{\partial M}$, such that the following axioms
hold:

\begin{itemize}
\item[A1)]  the assignment is functorial with respect to diffeomorphisms of $%
\Sigma $ and $M$,

\item[A2)]  $V_{-\Sigma }=V_{\Sigma }^{*}$, where $-\Sigma $ denotes $\Sigma 
$ with the opposite orientation, and $V_{\Sigma }^{*}$ is the dual vector
space of $V_{\Sigma }$,

\item[A3)] 
\begin{itemize}
\item[a)]  $V_{\Sigma _{1}\sqcup \Sigma _{2}}=V_{\Sigma _{1}}\otimes
V_{\Sigma _{2}}$, where $\sqcup $ denotes the disjoint union,

\item[b)]  $Z_{M_{1}\sqcup M_{2}}=Z_{M_{1}}\otimes Z_{M_{2}}$,

\item[c)]  (gluing axiom) $Z_{M_{1}\sqcup _{\Sigma
}M_{2}}=<Z_{M_{1}},Z_{M_{2}}>$, where $M_1\sqcup_\Sigma M_2$ means $M_1$ glued
to $M_2$ along $\Sigma$  (see figure) and 
$<.,.>$ is given by evaluation of $%
V_{\Sigma }^{*}$ on $V_{\Sigma }$,
\end{itemize}
\item[A4)]  a number of non-triviality conditions.
\end{itemize}

\begin{figure}[h]
\centerline{\psfig{figure=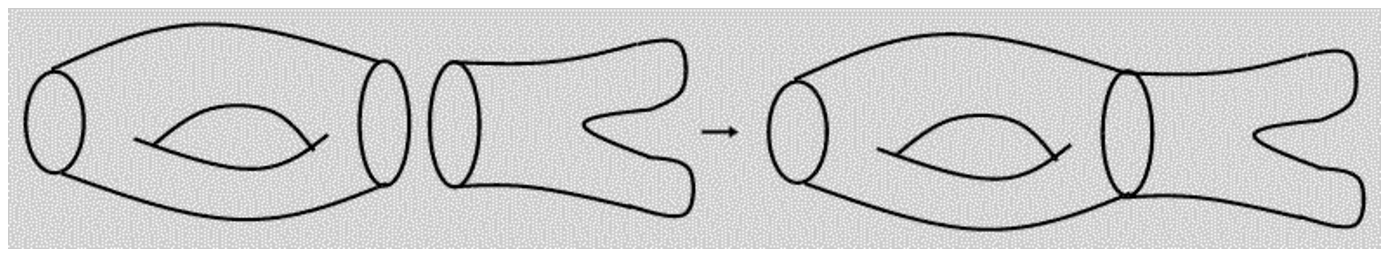,height=2.8578cm,width=16.0617cm}}
\end{figure}

The first axiom means, in particular, that when $M$ and $M^{\prime }$ are
diffeomorphic, there is a linear isomorphism $V_{\partial M}\rightarrow
V_{\partial M^{\prime }}$ which takes $Z_{M}$ to $Z_{M^{\prime }}$. One of
the non-triviality conditions is $V_{\emptyset }=\Bbb{K}$, which is related
to axiom A3a): $\emptyset \sqcup \emptyset =\emptyset \Rightarrow
V_{\emptyset }\otimes V_{\emptyset }=V_{\emptyset }$, hence $V_{\emptyset }$
can only be $\Bbb{K}$ or the trivial vector space $\{0\}$. Thus, in
particular, a TQFT assigns to every closed $(d+1)$ manifold $M$ a numerical
invariant $Z_{M}\in \Bbb{K}$.

The all-important gluing axiom may be reformulated by noting that, if $%
\partial M=\Sigma _{1}\sqcup (-\Sigma _{2})$, then $Z_{M}\in V_{\Sigma
_{1}}\otimes V_{\Sigma _{2}}^{*}\cong \mathrm{Lin}(V_{\Sigma _{2}},V_{\Sigma
_{1}})$. Thus $Z_{M}$ may be regarded as a linear map, and with $\partial
M_{1}=\Sigma _{1}\sqcup (-\Sigma )$ and $\partial M_{2}=\Sigma \sqcup
(-\Sigma _{2})$, the gluing axiom 3c) reads: 
\[
Z_{M_{1}\sqcup _{\Sigma }M_{2}}=Z_{M_{1}}\circ Z_{M_{2}}. 
\]
We remark that $\partial M$ can always be regarded as the disjoint union of
two boundary components, since either or both components can be empty.

At this stage it is probably helpful to give a simple example. When $d=0$,
any non-empty oriented $(d+1)$-dimensional manifold with boundary $M$ is
isomorphic to either the oriented interval $I$ or the oriented circle $S^{1}$%
, or to a disjoint union of copies of these. Suppose $\Bbb{K}$ is the field
of complex numbers, $\Bbb{C}$. To specify the TQFT, all we have to do is
choose $V_{\bullet +}$, the vector space assigned to the positively oriented
point $\bullet +$, and we choose it to be $\Bbb{C}^{2}$. All the remaining
assignments follow from the axioms. From A2), $V_{\bullet -}=(\Bbb{C}%
^{2})^{*}$, and the disjoint union of positively and negatively oriented
points maps to the corresponding tensor product of copies of $\Bbb{C}^{2}$
and its dual. Interpreting $Z_{I}$ as a linear map from $\Bbb{C}^{2}$ to $%
\Bbb{C}^{2}$, the topological diffeomorphism between two intervals glued
together at one end and a single interval, gives rise to the equation $%
Z_{I}\circ Z_{I}=Z_{I}$, and thus $Z_{I}$ is a projection. By another of the
non-triviality axioms, the projection $Z_{I}$ is actually taken to be
surjective, i.e., $Z_{I}=\mathrm{id}_{\Bbb{C}^{2}}$.

All that remains is to identify $Z_{S^{1}}$, which is a linear map from $%
V_{\emptyset }=\Bbb{K}$ to itself, and thus given by the number $%
Z_{S^{1}}(1)\in \Bbb{K}$. To this end, we look at the interval and $Z_{I}$
in two new ways. First, if we regard $\partial I$ as $-({\bullet +}\sqcup {%
\bullet -})\sqcup \emptyset $, by folding the interval into an upturned U
shape, oriented clockwise, and ``reading'' it from bottom to top (see
figure), 

\begin{figure}[h]
\centerline{\psfig{figure=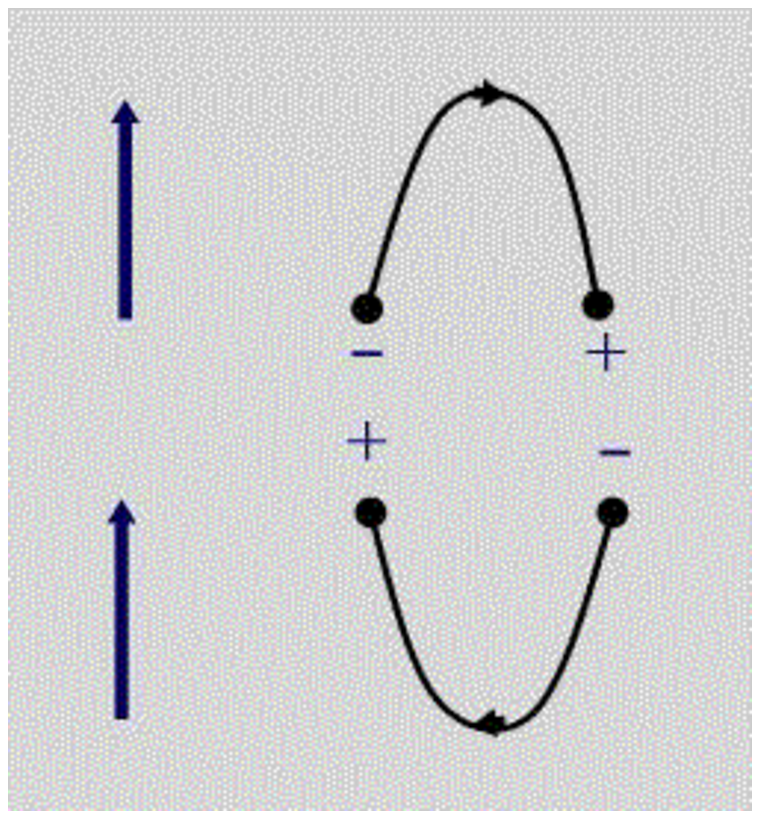,height=130.5625pt,width=120.9375pt}}
\end{figure}


\noindent we get $Z_{\cap +}:V_{\bullet +}\otimes V_{\bullet
-}\rightarrow \Bbb{C}$, where we have denoted the folded interval with its
clockwise orientation $\cap +$. In the same manner we may fold the interval
the other way into a clockwise-oriented U shape (see figure above) and
regard its boundary as $\emptyset \sqcup ({\bullet +}\sqcup {\bullet -})$.
In this way we get $Z_{\cup +}:\Bbb{C}\rightarrow V_{\bullet +}\otimes
V_{\bullet -}$, where we have denoted the folded interval with its clockwise
orientation ${\cup +}$. To describe the various linear transformations, we
introduce a basis $\{e_{1},e_{2}\}$ in $V_{\bullet +}=\Bbb{C}^{2}$, and let $%
\{e_{1}^{*},e_{2}^{*}\}$ denote its dual basis. The element of $V_{\bullet
+}\otimes V_{\bullet +}^{*}\cong \mathrm{Lin}(V_{\bullet +},V_{\bullet +})$
giving rise to $Z_{I}=\mathrm{id}_{\Bbb{C}^{2}}$ is $e_{1}\otimes
e_{1}^{*}+e_{2}\otimes e_{2}^{*}$. This same element, regarded as belonging
to $(V_{\bullet +}\otimes V_{\bullet +}^{*})\otimes \Bbb{C}$, leads to $%
Z_{\cup +}$ being given by $1\mapsto e_{1}\otimes e_{1}^{*}+e_{2}\otimes
e_{2}^{*}$. Regarding $e_{1}^{*}\otimes e_{1}+e_{2}^{*}\otimes e_{2}$ as
belonging to $\Bbb{C}\otimes (V_{\bullet +}\otimes V_{\bullet
+}^{*})^{*}\cong \mathrm{Lin}(V_{\bullet +}\otimes V_{\bullet +}^{*},\mathbf{%
\ }\Bbb{C})$ leads to $Z_{\cap +}$ being given by $e_{i}\otimes
e_{j}^{*}\mapsto \delta _{ij}$, for $i,j,=1,2$. Now we may consider the
clockwise oriented circle $S^{1}$ as the result of gluing a clockwise
oriented inverted U and a clockwise oriented U (see figure).

\begin{figure}[h]
\centerline{\psfig{figure=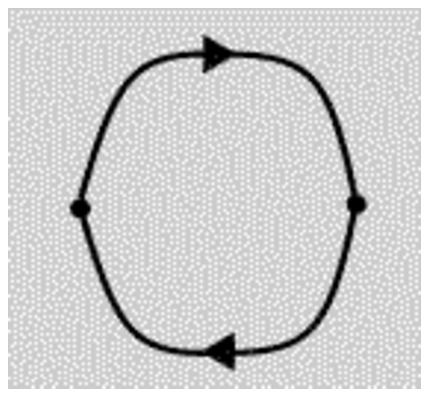,height=1.0914in,width=0.9513in}}
\end{figure}

\noindent Thus we calculate $Z_{S^{1}}(1)=Z_{\cap +}\circ Z_{\cup
+}(1)=Z_{\cap +}(e_{1}\otimes e_{1}^{*}+e_{2}\otimes e_{2}^{*})=2$. This
result completes the description of the TQFT, since any disjoint union of
intervals and circles maps to the appropriate tensor product of $Z_{I}$ and $%
Z_{S^{1}}$.

The preceding example also hints at yet another frequently-employed
reformulation of the Atiyah axioms for TQFT. This version may be regarded as
a response to the following, seemingly bizarre, question: if $M$ is mapped
to $Z_{M}$, a (linear) function, why not make $M$ itself into a function? In
fact, this is achieved by introducing the cobordism category, $(d+1)$-$%
\mathbf{Cobord}$, whose objects are oriented $d$-dimensional manifolds and
whose morphisms are equivalence classes of $(d+1)$-dimensional manifolds
with boundary. More precisely, if $\Sigma _{1}$ and $\Sigma _{2}$ are two
objects, then the set of all morphisms from $\Sigma _{2}$ to $\Sigma _{1}$
is $\limfunc{Mor}(\Sigma _{2},\Sigma _{1})=\left\{ M|\,\partial M=\Sigma
_{1}\sqcup (-\Sigma _{2})\right\} /\mathit{Diff}$, where $\mathit{Diff}$
stands for diffeomorphisms which restrict to the identity on the boundary of 
$M$. Composition of $[M_{1}]\in \mathrm{Mor}(\Sigma ,\Sigma _{1})$ and $%
[M_{2}]\in \mathrm{Mor}(\Sigma _{2},\Sigma )$ is given by $[M_{1}]\circ
[M_{2}]=[M_{1}\sqcup _{\Sigma }M_{2}]$, where the square brackets denote the
equivalence class. This is associative, due to the identification of
diffeomorphic manifolds, and furthermore, every object $\Sigma $ has an
identity morphism, namely $[\Sigma \times I]$.

In this framework, a TQFT is a functor from the cobordism category $(d+1)$-$%
\mathbf{Cobord}$ to the category of finite-dimensional vector spaces (over $%
\Bbb{K}$) $\mathbf{Vect}$, described by $\Sigma \mapsto V_{\Sigma }$ on
objects and $(\Sigma _{2}\stackrel{[M]}{\rightarrow }\Sigma _{1})\mapsto
(V_{\Sigma _{2}}\stackrel{Z_{[M]}}{\rightarrow }V_{\Sigma _{1}})$ on
morphisms, which preserves the products or ``monoidal structures'' ($\sqcup $
and $\otimes $ respectively), the unit objects for these products ($%
\emptyset $ and $\Bbb{K}$ respectively), and the ``involutions'' ($\Sigma
\mapsto -\Sigma $ and $V\mapsto V^{*}$ respectively).

This ``cobordism definition'' of TQFT is very elegant, but slightly unusual
from the viewpoint of conventional algebraic topology in that the cobordism
category is somewhat uncanonical as a category. We will return to discussing
this point in section \ref{altapp}. It is very appropriate for describing
TQFT-like constructions involving embedded manifolds, like braids and
tangles, such as Turaev's operator-valued invariant for tangles \cite{Turaev}%
. See \cite{Picken} for a discussion by one of the authors of the cobordism
approach for these embedded TQFT's.

\section{A comparison with homology \label{comphol}}

At this stage it is appropriate to return to the theme of the title, and
inquire about the role of TQFT in the context of classical algebraic
topology. It is illuminating to examine some of the similarities and
differences between Atiyah's axioms for TQFT and the Eilenberg-Steenrod
axioms for homology, as described in many textbooks on algebraic topology,
e.g. \cite{Rotman}.

\begin{center}
\begin{tabular}[t]{|c|p{3in}|p{3in}|}
\hline
& \multicolumn{1}{|c|}{\textbf{\ \rule[0.2in]{0in}{0in}Homology}} & 
\multicolumn{1}{|c|}{\textbf{\ TQFT }} \\ 
& \multicolumn{1}{|c|}{\rule[0.2in]{0in}{0in}\textbf{Eilenberg-Steenrod
axioms }} & \multicolumn{1}{|c|}{\textbf{Atiyah axioms }} \\ \hline
1 & \rule[0.2in]{0in}{0in}The topological objects are pairs of topological
spaces $(X,Y)$ with $Y\subset X$.\medskip & The topological objects are
pairs $(M,\Sigma )$ with $\Sigma =\partial M\subset M$. \\ \hline
2 & \rule[0.2in]{0in}{0in}A topological object $(X,Y)$ is sent to an abelian
group 
\[
H(X,Y)=\tbigoplus\nolimits_{n}H_{n}(X,Y)\text{.} 
\]
& A topological object $(M,\Sigma )$ is sent to a vector space (and a point
belonging to it) $(V_{\Sigma },Z_{M})$. \\ \hline
3 & \rule[0.2in]{0in}{0in}The theory is additive in the sense that $\sqcup $
maps to $\oplus $.\medskip & The theory is multiplicative in the sense
that $\sqcup $ maps to $\otimes $. \\ \hline
4 & \rule[0.2in]{0in}{0in}One of the key axioms is an excision property,
i.e. related to subtracting one space from another one (see figure
below).\medskip & One of the key axioms is a gluing property, i.e. related
to adding one space to another one (see figure below). \\ \hline
5 & \rule[0.2in]{0in}{0in}A single point `$\bullet $' maps to $H(\bullet
,\emptyset )=\Bbb{Z}$, the simplest non-trivial free abelian group. When $%
\bullet $ does not map to $\Bbb{Z}$, the theory is described as a
generalized homology theory. & The empty set $\emptyset $ maps to $%
V_{\emptyset }=\Bbb{K}$, the simplest non-trivial vector space over $\Bbb{K}$%
. For $d=0$, the single point does not in general map to $\Bbb{K}$, and for $%
d>0$ it does not even make sense to speak of $V_{\bullet }$. \\ \hline
6 & \rule[0.2in]{0in}{0in}The theory is not geared to any specific
dimension. There is a connecting homomorphism relating $H_{n}(-,-)$ to $%
H_{n+1}(-,-)$.\medskip & The theory is (usually) geared to a specific
dimension. There are only hints of a dimensional ladder, linking TQFT's for
different dimensions. \\ \hline
7 & \rule[0.2in]{0in}{0in}Homology is a functor and the topological
morphisms are canonical maps. & TQFT is a functor only in the cobordism
approach, and there the topological morphisms are equivalence classes of
manifolds. \\ \hline
8 & \rule[0.2in]{0in}{0in}The theory can be applied to various topological
categories.\medskip & The theory, as it stands, applies only to
(differentiable) manifolds. \\ \hline
\end{tabular}
\end{center}

\begin{figure}[h]
\centerline{\psfig{figure=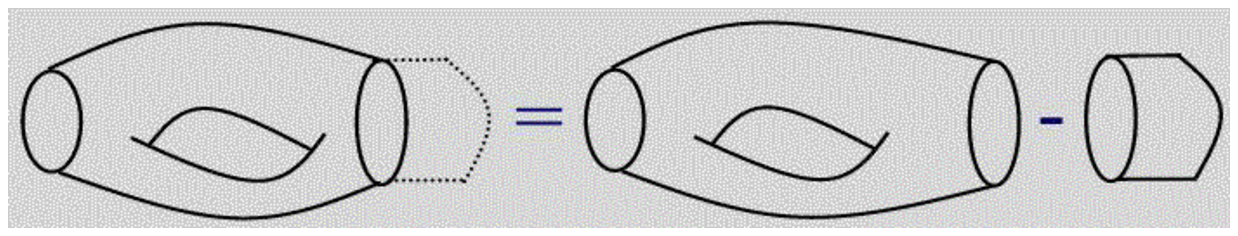,height=0.9651in,width=4.9121in}}
\end{figure}

\newpage

\begin{figure}[h]
\centerline{\psfig{figure=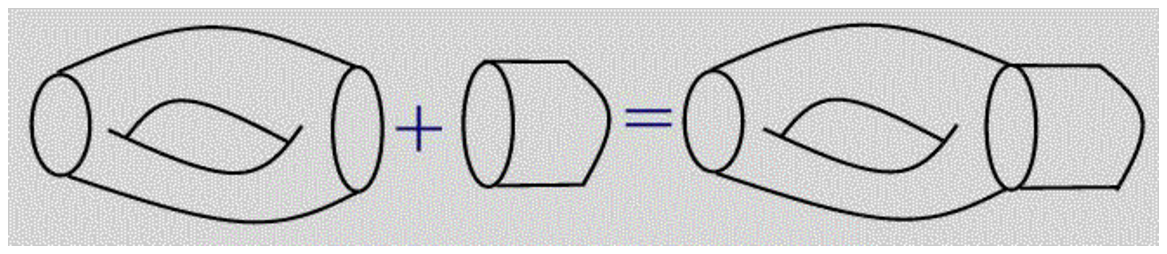,height=0.9305in,width=5.6637in}}
\end{figure}

We leave these comparisons for the reader to
contemplate, but cannot refrain from mentioning the observation, due to
Louis Crane, that TQFT behaves in some ways like an exponentiated version of
homology, bearing in mind especially properties 3 and 5.

The restriction to differentiable manifolds in the Atiyah axioms should not
be taken too literally, since, as Atiyah himself stresses, the axioms are
meant to be minimal and can be extended in a variety of ways. Indeed there
are a number of important examples of TQFTs for other topological
categories. We will proceed to describe two of these.

First there is a class of TQFTs known as ``state sum models'', which are
defined combinatorially using triangulated manifolds. One starts by
assigning algebraic data to the simplices of the triangulation, subject to
some admissibility conditions. These algebraic data may come from a variety
of sources, e.g. groups, quantum groups, representations, categories,
subfactors and so on. From the data one calculates a numerical weight by
means of some rule, and then, given a triangulated manifold $M$, chosen to
be without boundary for simplicity, $Z_M$ is defined to be the sum over all
admissible data of the corresponding weights. If the data and rules for
assigning weights are suitably matched, $Z_M$ is independent of the
triangulation chosen.

A very simple example of a state sum model, which gives the flavour of the
construction, is defined for triangulated $2$-manifolds without boundary in
the following fashion (this example is based on a construction of Dijkgraaf
and Witten \cite{DijWitten}). For any $2$-simplex of the triangulation
``colour'' its oriented edges with elements of a fixed finite group $G$,
subject to the admissibility condition that the group elements corresponding
to the $1$-cycle around the boundary of the $2$-simplex multiply to $1$ (the
identity). Also, if the orientation of an edge is reversed the group element
assigned to it is replaced by its inverse (see figure).

\begin{figure}[h]
\centerline{\psfig{figure=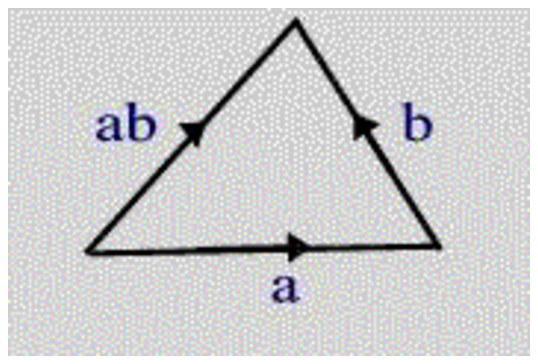,height=0.8043in,width=1.1978in}}
\end{figure}

\noindent Now, given a $2$-manifold $M$, we define: 
\[
Z_{M}=\sum_{\mathrm{colourings}}\left( \frac{1}{\#{G}}\right) ^{^{\#V}} 
\]
where $\#G$ is the number of elements of $G$ and $\#V$ is the number of
vertices of the triangulation. Triangulation invariance corresponds to
invariance under two local moves, the Pachner moves \cite{Pachner}, shown in
the figure below.

\begin{figure}[h]
\centerline{\psfig{figure=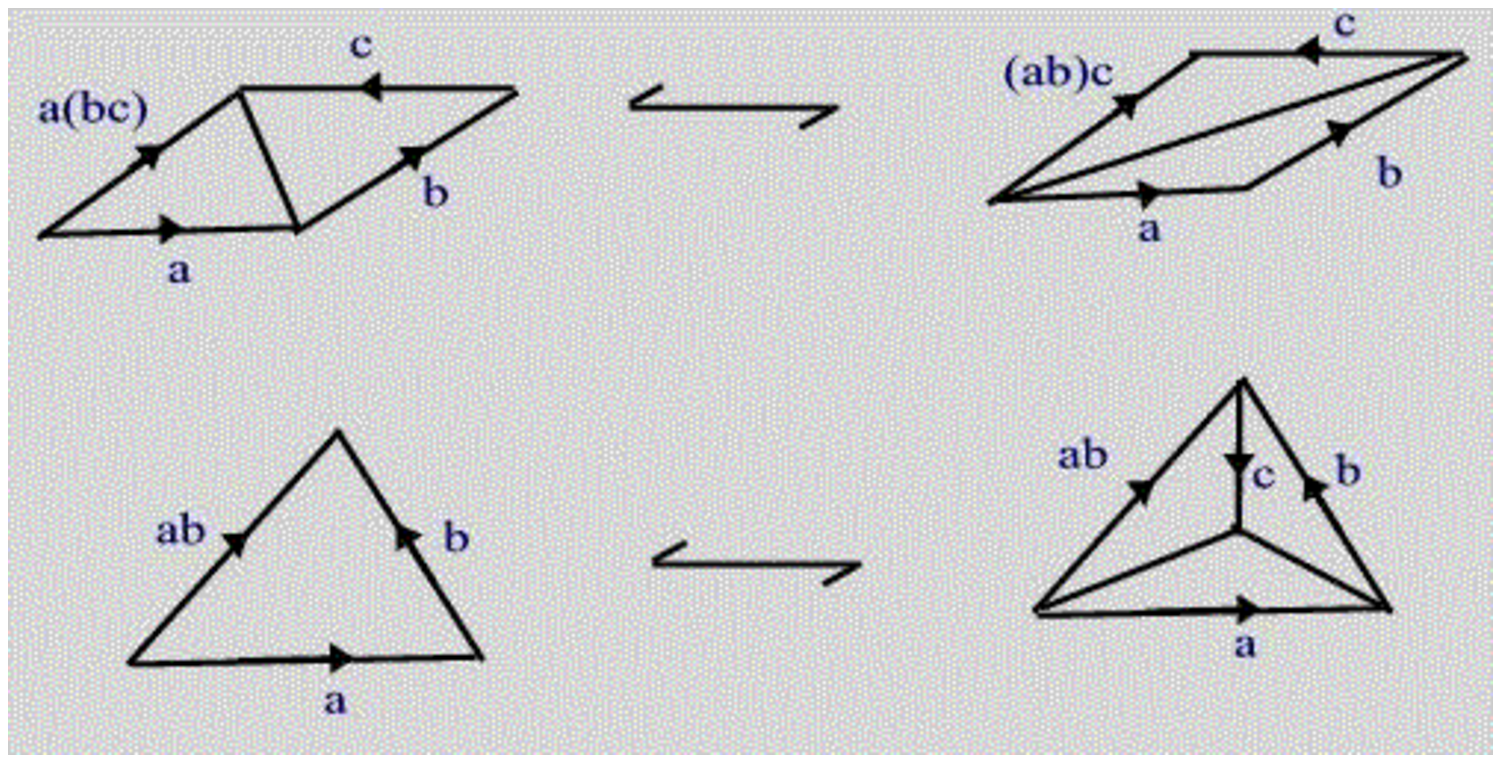,height=1.7037in,width=3.371in}}
\end{figure}

\noindent Invariance under these moves follows, for the first move, from the
associativity of group multiplication, and for the second move, from the
fact that the increase in the number of colourings by a factor of $\#G$, is
compensated by the $\frac{1}{\#{G}}$ factor coming from the extra vertex.

The second example of a TQFT not involving differentiable manifolds, is due
to Quinn \cite{Quinn}. This TQFT is defined for a very general class of
topological spaces, namely finite $CW$-complexes, which need not even be
(topological) manifolds. Let $M$ be a finite $CW$-complex, and let $\Sigma
_{1}^{i}$, $\Sigma _{2}^{o}$ be disjoint $CW$-subcomplexes of $M$, labelled $%
i$ for ``in'' and $o$ for ``out''. In the cobordism picture of TQFTs we can
regard $M$ (or rather an equivalence class of $M$s -- for notational
simplicity we will just write $M$) as a morphism from $\Sigma _{1}^{i}$ to $%
\Sigma _{2}^{o}$. If the ``out'' of $M_{2}$ equals the ``in'' of $M_{1}$,
both equal to $\Sigma $ say, we can glue or compose: 
\[
(\Sigma ^{i}\stackrel{M_{1}}{\longrightarrow }\Sigma _{1}^{o})\circ (\Sigma
_{2}^{i}\stackrel{M_{2}}{\longrightarrow }\Sigma ^{o})=(\Sigma _{2}^{i}%
\stackrel{M_{1}\sqcup _{\Sigma }M_{2}}{\longrightarrow }\Sigma _{1}^{o}), 
\]
where $M_{1}\sqcup _{\Sigma }M_{2}$ denotes $M_{1}$ and $M_{2}$ glued
together along their common subcomplex $\Sigma $.

Given this setup, Quinn defines a TQFT as follows. Choosing $\Bbb{K}=\Bbb{C}$%
, the complex numbers, we set $V_{\Sigma }=\Bbb{C}$ for every $\Sigma $,
i.e. this TQFT does not distinguish anything at the level of objects of the
cobordism category. (In fact, Quinn chooses a commutative ring $R$ instead
of $\mathbf{\ }\Bbb{C}$, but this makes little difference to the example.)
At the level of morphisms, however, the assignments are as follows: 
\[
(\Sigma ^{o}\stackrel{M}{\longrightarrow }\Sigma ^{i})\mapsto (\Bbb{C}%
\stackrel{Z_{M}}{\longrightarrow }\Bbb{C}),\quad Z_{M}(c)=e^{i\alpha \chi
(M,\Sigma ^{i})}(c) 
\]
where $\chi (M,\Sigma ^{i})=\sum_{n=0}^{\dim M}(-1)^{n}\limfunc{rank}%
(H_{n}(M,\Sigma ^{i}))$ is the relative homology version of the Euler
number, and $e^{i\alpha }$ is a fixed element of $\Bbb{C}$ of modulus $1$.
This assignment gives rise to a functor, and in particular it respects
composition: ($Z_{M_{1}\sqcup _{\Sigma }M_{2}}=Z_{M_{1}}\circ Z_{M_{2}}$),
because of the formula 
\[
\chi (M_{1}\sqcup _{\Sigma }M_{2},\Sigma _{2}^{i})=\chi (M_{1},\Sigma
^{i})+\chi (M_{2},\Sigma _{2}^{i}), 
\]
which, incidentally, may be proved by using excision. Another feature of
this example, relating to our comparison between TQFT and homology in the
previous section, is that it illustrates well the observation about the
exp(homology) structure of TQFT. We will be returning to a discussion of
this interesting example from a different perspective in section \ref
{quinnex}.

\section{Other definitions of TQFT \label{otherdef}}

The definition Atiyah gave in \cite{Atiyah} encapsulated the essential
ingredients common to a large class of TQFT models, whilst at the same time
restricting itself to a specific topological category, the category of
differentiable manifolds. As Atiyah himself states, the axioms allow for
numerous generalizations. Here we mention three other definitions which aim
to provide a more general framework for TQFT.

The first one, due to Quinn \cite{Quinn}, introduces the notion of ``domain
category'', being a category endowed with a collection of structures which
are abstractions of topological notions, such as boundary, cylinder, or
gluing. In particular the boundary of an object need not be the actual
boundary, since the object in question need not be a manifold, as in Quinn's
example described above. The axioms for a domain category are so general
that even purely algebraic examples, involving algebras over a commutative
ring, fit the definition.

Next, Turaev's definition, which appears in Chapter III of \cite{Turaev},
achieves generality in the topological category in a rather different way,
by introducing the abstract notion of ``space structure'', which encompasses
as special cases any kind of extra structure with which a topological space
can be endowed, such as a choice of orientation, a differentiable structure,
the structure of a CW complex, etc. Both Quinn and Turaev adopt the
``cobordism'' approach as described above, i.e. equivalence classes of $M$%
's are the morphisms of the topological cobordism category, gluing
corresponds to composition, and a TQFT is a functor from this cobordism
category to a suitable algebraic category.

The cobordism theme was taken a step further by Baez and Dolan in \cite
{BaezDolan}, when they started a programme to understand the subtle
relations between certain TQFT models for manifolds of different dimensions,
frequently referred to as the dimensional ladder. This programme is based on
higher-dimensional algebra, a generalization of the theory of categories and
functors to $n$-categories and $n$-functors, where for instance a $2$%
-category has not just objects and morphisms, but also $2$-morphisms, being,
roughly speaking, morphisms between morphisms. In this framework a TQFT
becomes an $n$-functor from the $n$-category of $n$-cobordisms to the $n$%
-category of $n$-Hilbert spaces. Since the definition of $n$-category is
itself rather elusive, this programme should still be described as being in
a state of development, but nevertheless reveals a fascinating view of
parallel developments in algebra and topology.

All three definitions are very much in the spirit of Atiyah's original
approach, which has thus proved to be rather influential. Also in practice
most authors studying a specific TQFT model base themselves on the cobordism
version of the Atiyah axioms.

\section{An alternative approach to TQFT \label{altapp}}

As discussed in the previous section, most authors adopt the cobordism
approach to TQFT, and thereby move away from the kind of framework which
might facilitate a comparison with conventional functors of algebraic
topology, such as homology and homotopy. Although the cobordism category is
indeed a category, its morphisms are not ``canonical maps'', whereas in the
topological categories used to define homology and homotopy theories they
are. More precisely, there is no forgetful functor from the cobordism
category to the category of topological spaces. One practical consequence is
that, in the cobordism framework, handling isomorphisms between $M$'s or $%
\Sigma $'s becomes somewhat delicate, since the role of morphisms in the
category has already been occupied by (equivalence classes of ) $M$'s. A
related observation is that composition of morphisms on the topological side
and the gluing operation are inextricably related in the cobordism approach.

Another point we would like to make is that, in the cobordism approach, the
gluing operation is inherently binary, i.e. necessarily involves gluing two
distinct $M$'s together. The possibility of gluing a single $M$ to itself
only enters at a later stage of development of the theory, for instance with
the following general result about gluing the ends of a ``generalized
cylinder'' $\Sigma \times I$, together, to make $\Sigma \times S^{1}$ (see
figure)

\begin{figure}[h]
\centerline{\psfig{figure=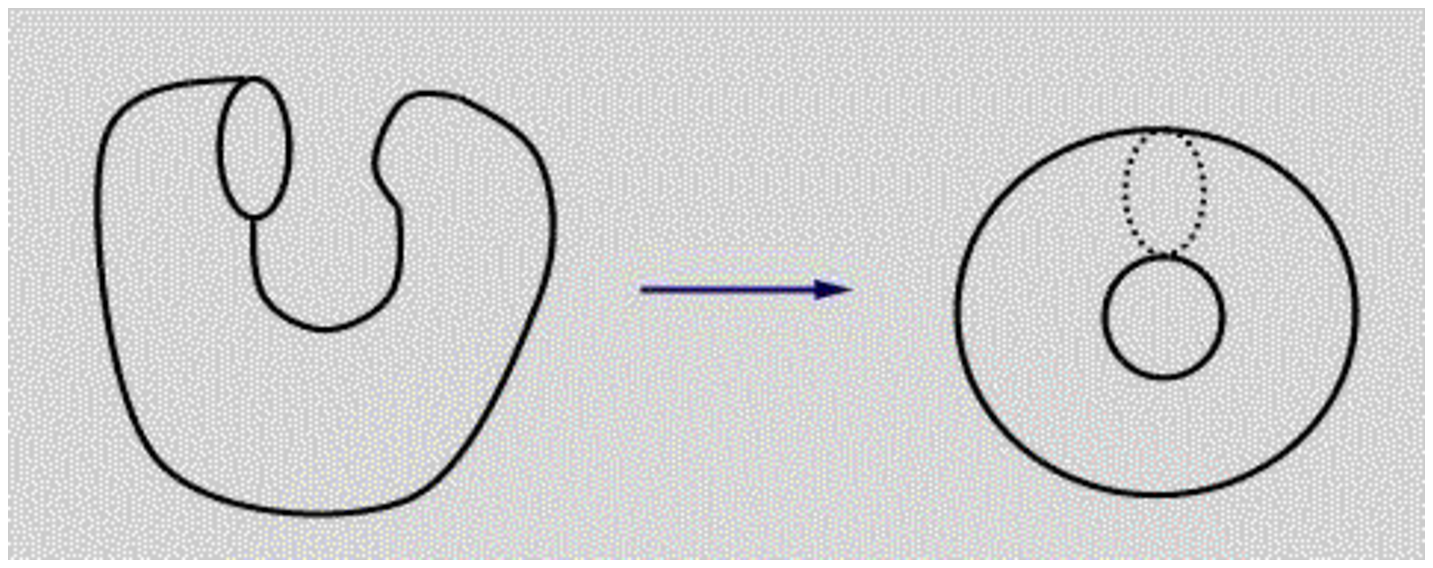,height=1.2635in,width=3.218in}}
\end{figure}

\noindent (here $ \Sigma $ is a $d$-dimensional manifold without boundary),
namely:  \[
Z_{\Sigma \times S^{1}}=\dim V_{\Sigma }. 
\]
We saw one instance of this TQFT theorem in section 2, in the example where $%
\Sigma $ was a single point. However, a more fundamental view of gluing is
as a \emph{unary} operation, which may of course be interpreted as binary
when the single $M$ actually consists of two separate parts.

These considerations led us to seek a different approach to defining TQFT
theories, with the following aims:

\begin{itemize}
\item  to incorporate a wide range of topological categories, as in the
Quinn and Turaev definitions discussed in the previous section, and
including cases of ``embedded topology'' like tangles;

\item  to formulate TQFT as a functor from a topological category, whose
morphisms are genuine maps in the above sense;

\item  to separate the roles of the composition and gluing operations;

\item  to incorporate gluing as a fundamentally unary operation from the
beginning;

\item  to provide a clean and efficient formulation for ease of calculation.
\end{itemize}

We will now proceed to describe the main features of this construction,
without going into technicalities. A detailed version is currently under
preparation. In a nutshell, the idea is as follows: starting with some
topological category, restrict the morphisms to be just isomorphisms and
so-called gluing morphisms, to be described shortly. A TQFT is then a
functor from this category to an algebraic category with suitably matching
structures.

The objects of the topological category $\mathcal{C}$ are \emph{pairs} of
the form $(M,\Sigma )$, where $M$ is an object of the aforementioned
starting category, denoted $\mathbf{M}$, e.g. the category of oriented $%
(d+1) $-dimensional manifolds with boundary, and $\Sigma $ is an object in a
related category, denoted $\mathbf{\Sigma }$, whose objects are subspaces of 
$M$'s with any additional structures that implies, e.g. the category of
oriented $d$-dimensional manifolds without boundary. Thus $\Sigma $ plays
the role of the boundary of $M$, although it need not be the actual boundary
in all cases (e.g. in Quinn's example, see sections \ref{comphol} and \ref
{quinnex}). Regarded as topological spaces, i.e. ignoring any structures,
the objects of $\mathbf{\Sigma }$ are taken to be finite disjoint unions of
connected and mutually separated components. Both $\mathbf{M}$ and $\mathbf{%
\Sigma }$ are monoidal categories with product the disjoint union $\sqcup $.

Furthermore there is a functor $I$ from $\mathbf{\Sigma}$ to itself, which
can be thought of as ``change of orientation'', and is such that $\Sigma$
and $I(\Sigma)$ are the same as topological spaces. Thus $I$ only acts on
the structures, not on the space itself.

Turning to the morphisms of the category $\mathcal{C}$, first we have
isomorphisms which, for a pair of objects $(M, \Sigma)$ and $(M^{\prime},
\Sigma^{\prime})$, are given by isomorphisms $f:M\rightarrow M^{\prime}$,
such that $f|_\Sigma $ is an isomorphism in the category $\mathbf{\Sigma}$
from $\Sigma$ to $\Sigma^{\prime}$. The only other type of morphisms we will
consider in $\mathcal{C}$ are the gluing morphisms which we will now define. 
\vspace{0.3cm}

\noindent \textbf{Definition} Let $(M,\Sigma )$ and $(M^{\prime },\Sigma
^{\prime })$ be two objects of $\mathcal{C}$ and suppose $\Sigma _{1}$ and $%
\Sigma _{2}$ are disjoint non-empty components of $\Sigma $, each being the
disjoint union of one or more connected components of $\Sigma $. A \emph{%
gluing morphism} from $(M,\Sigma )$ to $(M^{\prime },\Sigma ^{\prime })$ is
a pair $(f,\varphi )$, where $f:M\rightarrow M^{\prime }$ is a morphism of $%
\mathbf{M}$, and $\varphi :\Sigma _{1}\rightarrow I(\Sigma _{2})$ is an
isomorphism of $\mathbf{\Sigma }$, such that

\begin{itemize}
\item[1)]  $f$ is surjective,

\item[2)]  $f|_{M\setminus (\Sigma _{1}\sqcup \Sigma _{2})}$ is injective,

\item[3)]  $f|_{\Sigma \setminus (\Sigma _{1}\sqcup \Sigma _{2})}$ is an
isomorphism onto $\Sigma ^{\prime }$,

\item[4)]  for every $y\in f(\Sigma _{1})$ there is a unique pair $%
(x,\varphi (x))\in \Sigma _{1}\times \Sigma _{2}$ such that $f(x)=f(\varphi
(x))=y$,

\item[5)]  $f(\Sigma _{1})\cap f(M\setminus (\Sigma _{1}\sqcup \Sigma
_{2}))=\emptyset $,
\end{itemize}

\noindent where in conditions 1), 2), 4) and 5) $f$ and $\varphi$ refer to
the set-theoretic mappings underlying the corresponding morphisms. \vspace{%
0.3cm}

The intuitive content of the definition is that we are gluing two
``boundary'' components $\Sigma _{1}$ and $\Sigma _{2}$ of $M$ together
using the isomorphism $\varphi $, and the gluing morphism is from $M$
``before gluing'' to a copy of $M$ ``after gluing'' (see figure).

\begin{figure}[h]
\centerline{\psfig{figure=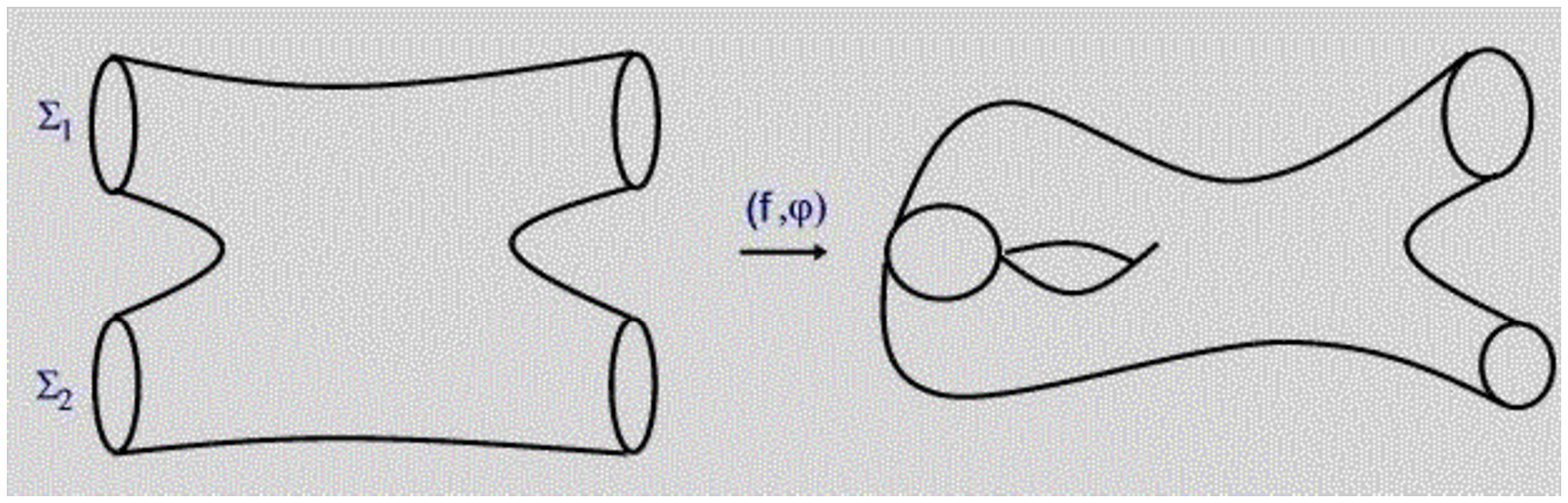,height=3.5cm,width=11cm}}
\end{figure}


From 2), $f$ is $1:1$ except on $\Sigma _{1}$ and $\Sigma _{2}$,
where, from 4) and 5), it maps $2:1$ onto their common image $f(\Sigma _{1})$%
, which, from 3) and 5), is disjoint from $\Sigma ^{\prime }$. Also, from
3), the unglued remainder of $\Sigma $ is isomorphic to $\Sigma ^{\prime }$.
In condition 4), as $\varphi $ is taken in the set-theoretic sense, the
distinction between $\Sigma _{2}$ and $I(\Sigma _{2})$ disappears, since
they are the same as sets.

The key property of gluing morphisms is that we can combine two of them to
get a new gluing morphism, defined to be their composition. \vspace{0.3cm}

\noindent \textbf{Theorem/Definition} Let $(f,\varphi )$ from $(M,\Sigma )$
to $(M^{\prime },\Sigma ^{\prime })$, and $(g,\psi )$ from $(M^{\prime
},\Sigma ^{\prime })$ to $(M^{\prime \prime },\Sigma ^{\prime \prime })$ be
gluing morphisms. Then $(g\circ f,\theta)$, where 
$\theta = \varphi \sqcup (I(f^{-1})\circ \psi \circ f)$, is a gluing
morphism from $(M,\Sigma )$ to $(M^{\prime \prime },\Sigma ^{\prime \prime })$, and is 
defined to be the composition of $(f,\varphi )$
and $(g,\psi )$. \vspace{0.3cm}

The intuitive content of the theorem is that gluing in two stages is
equivalent to gluing everything in one go.

\noindent \emph{Sketch of proof:} Let $\Sigma _{3}$ and $\Sigma _{4}$ be the
preimages under $f$ of $\Sigma _{1}^{\prime }$, $\Sigma _{2}^{\prime }$,
respectively, where $\psi :\Sigma _{1}^{\prime }\rightarrow I(\Sigma
_{2}^{\prime })$. By 3) and the fact that $\mathbf{\Sigma }$ is monoidal, $%
f|_{\Sigma _{3}}$ and $f|_{\Sigma _{4}}$ are isomorphisms onto $\Sigma
_{1}^{\prime }$, $\Sigma _{2}^{\prime }$, respectively. Then $\varphi \sqcup
(I(f^{-1})\circ \psi \circ f)$ is an isomorphism of $\mathbf{\Sigma }$,
since $\varphi $, $\psi $, $f|_{\Sigma _{3}}$ and $(f|_{\Sigma _{4}})^{-1}$
are all isomorphisms of $\mathbf{\Sigma }$, and $I$ is a functor. The
properties 1) to 5) are easily checked. \hfill $\Box $

It is straightforward to extend the above result to include the combination
of a gluing morphism and an isomorphism, or an isomorphism and a gluing
morphism, in both case giving rise to a new gluing morphism, defined as the
composition. Thus the isomorphisms and gluing morphisms taken together close
under composition. It is also easy to see that composition is associative
and that there is an identity morphism for each object $(M, \Sigma)$.
Furthermore the monoidal structures on $\mathbf{M}$ and $\mathbf{\Sigma}$
coming from the disjoint union give rise to a monoidal structure on $%
\mathcal{C}$. Thus we have: \vspace{0.3cm}

\noindent \textbf{Theorem} The above definitions yield a monoidal category $%
\mathcal{C}$, whose morphisms consist of isomorphisms and gluing morphisms.%
\vspace{0.3cm}

To get a TQFT functor we first need to choose an algebraic target category $%
\mathcal{D}$. This category is endowed with structures to match those of $%
\mathcal{C}$, and we describe them very briefly. The objects of $\mathcal{D}$
are pairs of the form $(V,x)$, where $V$ is a finite dimensional vector
space over a ground field $\Bbb{K}$, and $x$ is an element of $V$ (more
general choices are possible but we do not go into this here). The morphisms
from $(V,x)$ to $(W,y)$ in $\mathcal{D}$ are linear maps $f:V\rightarrow W$
such that $f(x)=y$. $\mathcal{D}$ is a monoidal category with product the
tensor product: $(V,x)\otimes (W,y)=(V\otimes W,x\otimes y)$. There is an
endofunctor $J$, corresponding to the endofunctor $I$ in $\mathcal{C}$,
which acts only on the vector space part. For instance, if $\Bbb{K}$ is $%
\Bbb{C}$, the complex numbers, $J$ could be the functor that maps $V$ to $%
\overline{V}$, being the space $V$ with the same addition as $V$ and the
conjugate scalar multiplication ($c\overline{\cdot }x:=\overline{c}\cdot x$).

A TQFT functor $Z$ is a functor from $\mathcal{C}$ to $\mathcal{D}$ of a
certain special form. To describe it we need the following assignments:

\begin{itemize}
\item[a)]  for each object $(M,\Sigma )$ of $\mathcal{C}$ an object $%
(V_{\Sigma },Z_{M})$ of $\mathcal{D}$, such that $V_{\Sigma }$ depends only
on $\Sigma $.

\item[b)]  for each isomorphism (of $\mathbf{\Sigma }$), $\Sigma _{1}%
\stackrel{\varphi }{\longrightarrow }\Sigma _{2}$, a linear isomorphism $%
V_{\Sigma _{1}}\stackrel{Z_{\varphi }}{\longrightarrow }V_{\Sigma _{2}}$,
with this assignment being functorial (mapping identity maps to identity
maps and compositions to compositions) and respecting the monoidal
structures and endofunctors $I$, $J$.

\item[c)]  for each $\Sigma $ of $\mathbf{\Sigma }$ a linear map (which we
call evaluation) $e_{V_{\Sigma }}:J(V_{\Sigma })\otimes V_{\Sigma
}\rightarrow \Bbb{K}$ satisfying various properties, including a
multiplicative property connecting $e_{V_{\Sigma \sqcup \Sigma ^{\prime }}}$
with $e_{V_{\Sigma }}$ and $e_{V_{\Sigma ^{\prime }}}$.
\end{itemize}

A TQFT functor $Z:\mathcal{C}\rightarrow \mathcal{D}$ is then given in terms
of these assignments by:

\begin{itemize}
\item[1)]  on objects: $(M,\Sigma )\mapsto (V_{\Sigma },Z_{M})$

\item[2)]  on isomorphisms: $((M,\Sigma )\stackrel{f}{\longrightarrow }%
(M^{\prime },\Sigma ^{\prime }))\mapsto ((V_{\Sigma },Z_{M})\stackrel{%
Z_{f|_{\Sigma }}}{\longrightarrow }(V_{\Sigma ^{\prime }},Z_{M^{\prime }}))$

\item[3)]  on gluing morphisms: 
\[
((M,\Sigma )\stackrel{(f,\varphi )}{\longrightarrow }(M^{\prime },\emptyset
))\mapsto ((V_{\Sigma },Z_{M})\stackrel{Z_{(f,\varphi )}}{\longrightarrow }(%
\Bbb{K},Z_{M^{\prime }}))\text{,} 
\]
where, for $\Sigma =\Sigma _{1}\sqcup \Sigma _{2}$ and $\varphi :\Sigma
_{1}\rightarrow I(\Sigma _{2})$, $Z_{(f,\varphi )}$ is given by 
\[
Z_{(f,\varphi )}=e_{V_{\Sigma _{2}}}\circ (Z_{\varphi }\otimes \mathrm{id}%
_{V_{\Sigma _{2}}})\text{,} 
\]
with a similar but somewhat more complicated formula for the case when $%
\Sigma ^{\prime }\neq \emptyset $, involving $e_{V_{\Sigma _{2}}}$, $%
Z_{\varphi }$ and $Z_{f|_{\Sigma \setminus (\Sigma _{1}\sqcup \Sigma _{2})}}$%
.
\end{itemize}

The main theorem is that the properties of the assignments a)-c) guarantee
that $Z$ is a functor from $\mathcal{C}$ to $\mathcal{D}$ respecting the
monoidal structures. Essentially $Z$ provides a representation of a severely
restricted subclass of morphisms of the starting category $\mathbf{M}$, but
a subset which includes the all-important class of gluing morphisms allowing
for topology changes. In the next section we will re-examine one of our
previous examples from this new perspective to illustrate how the definition
works in practice.

\section{Quinn's example revisited \label{quinnex}}

In this section we return to Quinn's example of a TQFT for finite $CW$%
-complexes, rephrased in terms of our definition. The objects of the
topological category $\mathcal{C}$ are pairs $(M,\Sigma )$, where $M$ is a
finite $CW$-complex and $\Sigma $ is a disjoint union of connected
subcomplexes of $M$, with each connected component labelled $i$ (in) or $o$
(out). We denote the union of the in (out) components of $\Sigma $ by $%
\Sigma ^{(i)}$ ($\Sigma ^{(o)}$). The functor $I$ acts on $\mathbf{\Sigma }$
by changing the labels from $i$ to $o$ and vice-versa, i.e. $I(\Sigma
^{i})=\Sigma ^{o}$ (or $I(\Sigma ^{o})=\Sigma ^{i}$), but $\Sigma ^{i}$ and $%
\Sigma ^{o}$ as $CW$-subcomplexes are the same. The morphisms of $\mathcal{C}
$ are isomorphisms and gluing morphisms, where isomorphisms are given by
isomorphisms in the category of finite $CW$-complexes $f:M\rightarrow
M^{\prime }$, such that $f|_{\Sigma }$ is an isomorphism from the subcomplex 
$\Sigma $ to $\Sigma ^{\prime }$ preserving the labels of each component,
and gluing morphisms are given by the general definition in the previous
section.

The objects of the algebraic category $\mathcal{D}$ are pairs of the form $%
(V,x)$ where $V$ is a vector space over the field of complex numbers $\Bbb{C}
$, and $x\in V$. The endofunctor $J$, corresponding to the topological
endofunctor $I$, is taken to be trivial, i.e. $J$ is the identity functor.
The morphisms of $\mathcal{D}$ are linear maps preserving the respective
elements.

A class of TQFT functors may now be specified as follows. For every $\Sigma $
we take $V_{\Sigma }=\Bbb{C}$. The element of $V_{\Sigma }$ corresponding to 
$M$ is: 
\[
Z_{M}=u^{c_{1}\chi (M)+c_{2}\chi (\Sigma ^{(i)})+c_{3}\chi (\Sigma ^{(o)})}%
\text{ ,} 
\]
where $u$ is some fixed nonzero element of $\Bbb{C}$, $\chi $ denotes the
Euler characteristic, and $c_{1},c_{2}$ and $c_{3}$ are fixed unknowns. The
assignment $\varphi \mapsto Z_{\varphi }$ is taken to be trivial, i.e. $%
Z_{\varphi }$ is the identity map on $\Bbb{C}$ for every isomorphism $%
\varphi $ of $\mathbf{\Sigma }$. Finally the evaluation $e_{V_{\Sigma
}}:J(V_{\Sigma })\otimes V_{\Sigma }=\Bbb{C}\otimes \Bbb{C}\rightarrow \Bbb{C%
}$ is given by $e_{V_{\Sigma }}(x\otimes y)=u^{c_{4}\chi (\Sigma )}xy$. The
multiplicative property mentioned in the general case, here corresponds to
the statement $\chi (\Sigma \sqcup \Sigma ^{\prime })=\chi (\Sigma )+\chi
(\Sigma ^{\prime })$.

From these assignments we can determine how $Z_{M}$ changes under
topological isomorphisms and gluing morphisms. For an isomorphism $%
f:(M,\Sigma )\rightarrow (M^{\prime },\Sigma ^{\prime })$ we get $%
Z_{M}=Z_{M^{\prime }}$, since $Z_{f|_{\Sigma }}$ is the identity. To study
the effect of gluing morphisms, let us start, in terms of the framework of
section \ref{comphol}, with the case of two separate $M$'s, $(M_{1},\Sigma
_{1}^{i}\sqcup \Sigma _{2}^{o})$ and $(M_{2},\Sigma _{2}^{i}\sqcup \Sigma
_{3}^{o})$, which are glued together along their common component $\Sigma
_{2}$ to give $(M_{1}\sqcup _{\Sigma _{2}}M_{2},\Sigma _{1}^{i}\sqcup \Sigma
_{3}^{o})$. In our terms this corresponds to a gluing morphism 
\[
(f,\varphi ):(M_{1},\Sigma _{1}^{i}\sqcup \Sigma _{2}^{o})\sqcup
(M_{2},\Sigma _{2}^{i}\sqcup \Sigma _{3}^{o})\rightarrow (M_{1}\sqcup
_{\Sigma _{2}}M_{2},\Sigma _{1}^{i}\sqcup \Sigma _{3}^{o}) 
\]
where $\varphi $ is the identity on $\Sigma _{2}$. Now $Z_{(f,\varphi )}$
acts here by multiplication by $u^{c_{4}\chi (\Sigma _{2})}$, due to our
choice of evaluation, so in the equation 
\[
Z_{(f,\varphi )}(Z_{M_{1}}\otimes Z_{M_{2}})=Z_{M_{1}\sqcup _{\Sigma
_{2}}M_{2}}\text{,} 
\]
we have the exponent of $u$ on the left hand side given by: 
\[
c_{4}\chi (\Sigma _{2}^{i})+c_{1}\chi (M_{1})+c_{2}\chi (\Sigma
_{1}^{i})+c_{3}\chi (\Sigma _{2}^{o})+c_{1}\chi (M_{2})+c_{2}\chi (\Sigma
_{2}^{i})+c_{3}\chi (\Sigma _{3}^{o}) 
\]
and on the right hand side by: 
\[
c_{1}\chi (M_{1}\sqcup _{\Sigma _{2}}M_{2})+c_{2}\chi (\Sigma
_{1}^{i})+c_{3}\chi (\Sigma _{3}^{o}). 
\]
Due to the formula $\chi (M_{1}\sqcup _{\Sigma _{2}}M_{2})=\chi (M_{1})+\chi
(M_{2})-\chi (\Sigma _{2})$, which already appeared in a slightly different
form in section \ref{comphol}, the equation is equivalent to $%
(c_{1}+c_{2}+c_{3}+c_{4})\chi (\Sigma _{2})=0$, and since $\chi (\Sigma
_{2}) $ is not necessarily zero, we get the following constraint on the
constants: 
\[
c_{1}+c_{2}+c_{3}+c_{4}=0\text{.} 
\]

Quinn in his original discussion \cite{Quinn} gave two examples of TQFTs,
the one we described in section \ref{comphol}, which he called the Euler
theory, and a modified example called the skew Euler theory. In terms of our
approach these correspond to two special cases of the above constraint: 
\[
\begin{tabular}{cccl}
$c_{1}=-c_{2}=1$ & and & $c_{3}=c_{4}=0$ & (Euler theory), \\ 
$c_{1}=-c_{3}=1$ & and & $c_{2}=c_{4}=0$ & (skew Euler theory).
\end{tabular}
\]
Replacing $u$ by $u^{\prime }=u^{c}$ means that one of the unknowns can be
set to $1$ and the most natural choice is to set $c_{1}=1$. So, from the
equation above we get 
\[
c_{2}+c_{3}+c_{4}=-1\text{.} 
\]
Apart from Quinn's two solutions there is a ``balanced solution'', with $%
c_{4}=0$ and $c_{2}=c_{3}=-1/2$, which is halfway between them in the sense
that the Euler characteristics of $\Sigma ^{(i)}$ and $\Sigma ^{(o)}$ appear
on an equal footing in the formula for $Z_{M}$, but there are many other
solutions, even with $c_{4}=0$. Any solution with $c_{4}\neq 0$ has the
property that the corresponding TQFT does not map the topological morphisms
trivially, i.e. does not map all morphisms to the identity map on $\Bbb{C}$%
.\smallskip

In our framework we can go one step further and consider the effect of
self-gluing. First we need a formula to replace the previous one for $\chi
(M_{1}\sqcup _{\Sigma _{2}}M_{2})$. Let $b_{n}(M)=\limfunc{rank}(H_{n}(M))$,
where $H_{n}(M)$ is finitely generated, since $M$ is a finite $CW$-complex.
We have 
\[
\chi (M)=\tsum\limits_{i=0}^{\infty }(-1)^{i}b_{i}(M). 
\]
Suppose we have disjoint subcomplexes $\Sigma _{1}$ and $\Sigma _{2}$ of $M$%
, and $\varphi :\Sigma _{1}\rightarrow \Sigma _{2}$ is an isomorphism. Let $M%
\stackrel{\nu }{\rightarrow }M_{\varphi }$ be the canonical map, where $%
M_{\varphi }$ is the identification space under the equivalence relation
generated by $x\sim \varphi (x)$ for any $x\in \Sigma _{1}$. The spaces $%
M\setminus (\Sigma _{1}\sqcup \Sigma _{2})$ and $M_{\varphi }\setminus \nu
(\Sigma _{1})$ are homeomorphic and thus $H_{n}(M\setminus (\Sigma
_{1}\sqcup \Sigma _{2}))$ is isomorphic to $H_{n}(M_{\varphi }\setminus \nu
(\Sigma _{1}))$. Now using excision we get $H_{n}(M,\Sigma _{1}\sqcup \Sigma
_{2})\cong H_{n}(M_{\varphi },\nu (\Sigma _{1}))$ and hence $%
b_{n}(M_{\varphi })=b_{n}(M)-b_{n}(\Sigma _{2})$. Thus, the Euler
characteristic formula for self-gluing is: 
\[
\chi (M_{\varphi })=\chi (M)-\chi (\Sigma _{2})\text{.} 
\]
Now in our TQFT approach, we have a gluing morphism $(f,\varphi ):(M,\Sigma
_{1}^{i}\sqcup \Sigma _{2}^{o})\rightarrow (M_{\varphi },\emptyset )$. The
equation $Z_{(f,\varphi )}(Z_{M})=Z_{M_{\varphi }}$ gives rise to the
equation for the exponents of $u$ on either side: 
\[
c_{4}\chi (\Sigma _{2})+\chi (M)+c_{2}\chi (\Sigma _{1})+c_{3}\chi (\Sigma
_{2})=\chi (M_{\varphi }) 
\]
and using $\chi (\Sigma _{1})=\chi (\Sigma _{2})$ and the above formula for $%
\chi (M_{\varphi })$, this corresponds to the same constraint as for mutual
gluings 
\[
c_{2}+c_{3}+c_{4}=-1\text{.} 
\]

It is straightforward to extend the previous discussion to cases where some
out components of $M_{1}$ and some in components of $M_{2}$ remain after
gluing, in the case of gluing two $M$s together, and some in and out
components remain after gluing, in the case of self-gluing of a single $M$%
.\smallskip

In conclusion, our approach to TQFT allows one to considerably increase the
class of Quinn-type examples and extend their range of application. It is
our hope that this approach will help clarify some issues in the general
theory and specific TQFT models, and will inspire new types of TQFT
construction. In future work, apart from giving full details of the
definition, we intend to develop other examples, including ones involving
embedded topology, like curves in manifolds, or some geometrical features.
It is our belief that TQFT is a very profound structure offering a wide
range of potential applications still to be explored.

\newpage

\bibliographystyle{abbrv}
\bibliography{art,artTQFT,LivrosAlgTop,LivrosCatHomAlg}

\end{document}